\documentclass[11pt]{article}
\usepackage{latexsym}
\usepackage{amsfonts}
\makeatletter
\def\abstract{\small\quotation{\hskip-\parindent\sc Abstract.}}

\def\classification{\@ifnextchar [{\@xfootnotenext}%
   {\begingroup\let\protect\noexpand
      \xdef\@thefnmark{}\endgroup
    \@footnotetext}}

\title {}

\begin{document}

\classification {{\it 1991 Mathematics
Subject Classification:} Primary 20F28,  secondary 20E05, 20F36.} 

\begin{center}
{\bf \large   FIXED  POINTS  OF  ENDOMORPHISMS  OF  A  FREE 
METABELIAN  GROUP}  
\bigskip

{\bf    Vladimir Shpilrain}

\end{center} 
\medskip

\begin{abstract}
\noindent  We consider IA-endomorphisms  (i.e.,  Identical  in  
Abelianization) of a free metabelian group of finite rank,  and 
give  a  matrix 
 characterization of their fixed points which  is   similar  to  
   (yet different from) the 
 well-known characterization of eigenvectors of a linear  
operator  in  a  vector space.    We  then use our matrix
 characterization to elaborate several properties  of 
 the fixed point groups of metabelian  endomorphisms.   In
particular, we show that the rank of  the  fixed  point group of an 
 IA-endomorphism of the free metabelian group of rank $n \ge 2$ 
 can be either  equal to 0, 1, or greater than 
$(n-1)$ (in particular, it can be  infinite). 
We also point out     a connection between these properties of
metabelian  IA-endomorphisms  and some properties of the Gassner
representation  of pure braid groups. 

\end{abstract}
\date{}

\bigskip

\noindent {\bf 1. Introduction }
\bigskip

 \indent       Let  $F = F_n$   be the free group of a finite rank 
$n \ge 2$  with a  set  
$X = {\{}x_i {\}}, 1 \le i \le n$,  of free generators.  Then, let  
$F' =  [F,  F]$   be 
 the commutator subgroup of the group  $F$.  The group  $M = M_n  =
F/F'' $  is  called a free metabelian group. 
\smallskip 

      Let  $\varphi$  be an endomorphism of the group  $F$  given
by   $\varphi(x_k) = y_k , 1 \le k \le n$. 
 It also induces an endomorphism of  the  group   $M$   in  the 
 natural way.  We denote elements of  a  free  group  and  their  
 natural 
 images in a free metabelian group by the same letters when there  
 is  no   ambiguity.  The same applies to endomorphisms. 
\smallskip 

      Birman \cite {Bi1} and Bachmuth \cite {Ba} have used 
matrices  over  group  
 rings  to  study  endomorphisms  of  a  free  and  a  free  
metabelian   group,  respectively. 
\smallskip 

      Birman \cite {Bi1} has given a  matrix  characterization  
 of  automorphisms  of a free group among arbitrary  endomorphisms 
(the  ``inverse  function  theorem") as follows.  Define the  matrix 
$J_{\varphi}  = (d_j (y_i))_{1 \le i,j \le n}$        (the 
 ``Jacobian matrix" of $\varphi$), where  $d_j$   denotes partial 
Fox
 derivation (with respect to $x_j$) 
in the  free 
 group ring  $ZF$  (see \cite {Fox}).  Then  $\varphi$  is an 
 automorphism  if and only  if  the matrix  $J_{\varphi}$   is
invertible. 
\smallskip 

      Bachmuth \cite {Ba} has obtained an inverse function theorem 
of  the  same 
 kind on replacing the  Jacobian  matrix   $J_{\varphi}$    by  
 its  image $~J_{\varphi}^{a}~$ over  the  abelianized group ring  
$Z(F/F')$.
 \smallskip 

      Recently, matrix methods have been used by a number of  
 authors  to  produce new interesting  results  on  endomorphisms 
of  free  and  free  metabelian groups.  Umirbaev \cite {U} has 
generalized  Birman's  result  to  partial generating systems 
(so-called  {\it primitive  systems})  of  a  free 
 group.  In \cite {GGN} and \cite {R},  similar  results  are 
obtained  for  primitive 
 systems of a free  metabelian  group.   In  \cite {Sh}, 
{\it monomorphisms}   (i.e.,  injective endomorphisms) of  a  free 
group  are  characterized  by  the  property of the Jacobian matrix
to have left independent rows. 
\smallskip 

      All these results show a remarkable parallelism between the 
 theory of  endomorphisms of a free (or a free metabelian) group and
 the  theory  of  linear operators in a vector space.  In this 
 paper,  we  obtain a matrix characterization of   
 IA-endomorphisms  with non-trivial fixed points (``eigenvectors")
which, although is similar to the corresponding well-known 
 characterization in linear algebra, also reveals a subtle
 difference.  
\smallskip 

 IA-endomorphism  of  a  group  $G$  is an endomorphism which is
 Identical  in  the  Abelianization  
 $G/G'$  of the group  $G$. As usual, by the  $rank$ of a matrix
$A$  over a commutative ring $S$ we mean the maximal number of rows
of  $A$ independent over $S$.

\medskip

\noindent {\bf  Theorem 1.1.} Let  $\varphi$  be an IA-endomorphism
of the free metabelian  group   $M = M_n$  given by  $\varphi(x_i) = 
y_i$,  $ ~1 \le i \le n,~$  and let 
$~J_{\varphi}^{a}~$   be   the  corresponding  abelianized  Jacobian 
matrix.   Then:
\smallskip 

\noindent {\bf (i)}  $ det(J_{\varphi}^{a}  - I) = 0$,   where  $I$ 
is the  $n$x$n$  identity matrix. 
\smallskip 

\noindent {\bf (ii)}  If $~rank(J_{\varphi}^{a}  - I)  \le n - 2$, 
then $\varphi$  has a non-trivial fixed point inside
the commutator subgroup $M'$.  
 \smallskip 

\noindent {\bf (iii)} If $~rank(J_{\varphi}^{a}  - I) = n - 1$, 
then $\varphi$ has a non-trivial fixed point inside
the commutator subgroup $M'$ if and  only if $[x_1, x_2]^{u_1}\cdot 
[x_2, x_3]^{u_2} \cdot ... \cdot [ x_{n - 1}, x_n]^{u_{n - 1}} = 
[y_1, y_2]^{u_1}\cdot
[y_2, y_3]^{u_2} \cdot ...  \cdot [ y_{n - 1}, y_n]^{u_{n - 1}}~$ 
for  some  elements  $u_k \in Z(M/M')$, not all of them zero. 

\medskip

 The situation is  most subtle when 
 $~rank(J_{\varphi}^{a}  - I) = n - 1$.    In this  case,       
             anything   
 can happen!  For example, if $\varphi$ is 
 an IA-endomorphism of the   group   $M_2$ given by $\varphi(x_1) = 
x_1 s; ~\varphi(x_2) = x_2 s^{-1}$, then $\varphi$ has no non-trivial 
fixed points for any $s \in M'_2, ~s \ne 1~$ (see Proposition 3.2).
On the other  hand, if $\varphi$ is an inner automorphism induced by
an element  $g \not\in M'_2$, then $\varphi$ has a non-trivial fixed
point (outside $M'_2$). In both situations, $~rank(J_{\varphi}^{a}  -
I) = 1 = 2 - 1$.
\smallskip 

Although I was not able to distinguish these  possibilities  
``by matrix means", an algorithm for detecting fixed points does 
exist:  
\medskip 

\noindent {\bf   Theorem  1.2.} There is an algorithm for
 detecting the presence of non-trivial 
fixed points of an arbitrary 
IA-endomorphism of a free metabelian  group $M_n$. Also, there 
is an algorithm for detecting the presence of those 
 non-trivial fixed points that belong to the commutator subgroup
$M'_n$. In both cases, if  an algorithm reveals the presence of
non-trivial  fixed points, then at least one of them can be actually
found. 
 \medskip 

In  case  of  a  free  group,    these  questions  remain  open 
even  when  restricted to
 detecting fixed points of {\it automorphisms}. 

\smallskip 

When all  non-trivial fixed points of $\varphi$ are 
 inside the commutator subgroup $M'$ (which is the case, for
example, when $\varphi$ is the conjugation by a non-trivial 
 element from $M'$), the fixed point group $Fix ~\varphi~$ has {\it
infinite rank}: 
 \medskip

\noindent {\bf   Proposition 1.3.}  Let  $\varphi$  be an 
IA-endomorphism  of     a free 
 metabelian group  $M$.  Then the intersection of  the  fixed  point 
group   $Fix ~\varphi~$  with  $M'$  is a  normal  subgroup  of 
$M$.   In  particular, if this  intersection is non-trivial, then
it is infinitely generated. 
\smallskip 

Furthermore, if $\varphi$ has the fixed  point group of
finite rank,  then  there is a ``generation gap": 
\medskip 

\noindent {\bf   Theorem  1.4.}  Let  $\varphi$  be an
IA-endomorphism  of     a free  metabelian group  $M_n$. If the
rank of  the  fixed  point group   $Fix ~\varphi~$ is finite, 
 then it is equal to 0, 1, or is greater than $~(n-1)~$. 
\medskip 

      Then, we have 
\smallskip 

\noindent {\bf   Proposition  1.5.}  For  an  arbitrary   $n \ge 3$,  
there   are   non-inner 
 (IA-)automorphisms of  the  free  metabelian  group    $M_n$    
that   have  infinitely generated fixed point group. 
\smallskip 

     Note that  every IA-automorphism of  the 
group   $M_2$    is inner \cite {Ba}.  Then, if  $u \in M'_n$, the 
conjugation by  $u $ is an 
 automorphism  of   $M_n$  whose  fixed  point  group  coincides 
with   $M'_n$   and  therefore  is 
 infinitely  generated.   The  same  argument  proves  the  
existence  of {\it inner} automorphisms with infinitely generated
fixed point group for any  group  of the form  $F/[R, R]$  with 
$F/R$  infinite. 

 We should 
also  mention here an
 example of  a {\it finitely  presented } group
\cite  {Turner1} whose  automorphisms can have infinitely generated 
fixed point group.
 \medskip

     All this makes contrast to 
 the situation  in  a  free 
group:   fixed  point group of any free endomorphism has rank  $n =
rank ~F $  or  less.   This has been proved in \cite {BH} for
automorphisms,  and    in \cite  {Turner2} - for  arbitrary
endomorphisms.  \smallskip 

      In general, it is an interesting and important question -  how  
far  is  the  similarity  between  free  and  free  metabelian 
endomorphisms  extended?   The desire to ``lift" properties  of 
metabelian  endomorphisms   
 to  those of free endomorphisms leads 
 to several interesting 
questions  of which the following one seems particularly
attrac- tive:  
\smallskip 

\noindent {\bf   Problem 1.}  Suppose  $\varphi \in  Aut ~F$  is an
IA-automorphism with  a  non-trivial  fixed point.  Is it true that 
$\varphi$  has a fixed point outside   $[F, F]$ ?  
 \medskip
 
 Another interesting question is - how can our Theorem 1.1 be 
    strengthened when restricted to {\it 
automorphisms}; in particular, we ask: 
\medskip

\noindent {\bf   Problem 2.} Is it true that every IA-automorphism 
 of the group  $M_n$ has a  non-trivial  fixed point? 
  \medskip

 Informally speaking, endomorphisms with non-trivial fixed points 
 are rare. On the other hand, endomorphisms which are
automorphisms, are also rare (in contrast to the situation 
in linear algebra!).  Whether or not automorphisms with non-trivial 
fixed points  are rare, remains a mystery.
\medskip 

 Finally, we point out  a connection between the  properties of
metabelian  IA-endo- morphisms described above,  and some properties
of the Gassner representation  of pure braid groups.  For
defenitions of braids, braid groups, and for a background material,
we refer to  \cite {Bi2}. 

 Let $~\sigma \in P_n$ be a pure braid on $n$ strands, and $~\hat
{\sigma}~$  the corresponding link. Then, let {\small {\small
$\Gamma$}}$_{n-1} (\sigma)$  be the image of  $~\sigma~$ under the
reduced Gassner representation  of the pure braid group $P_n$, and
$A_{\hat \sigma}(t_1, ..., t_n)$  the Alexander polynomial of
$~\hat{\sigma}$. Then (see  [4, Theorem 3.11]): 
\begin{center}
  $A_{\hat \sigma}(t_1,..., t_n) = 0~$  if and only if  
 $~det(${\small {\small $\Gamma$}}$_{n-1} (\sigma)   - I) = 0$.
\end{center} 

 To explain how this is connected to the subject of the present
paper, we mention that the pure braid group $P_n$ is isomorphic 
to a subgroup of the group of IA-automorphisms of $F_n$, and 
the Gassner representation maps every IA-automorphism in this
subgroup onto its  abelianized Jacobian  matrix. Then,   
the  matrix {\small {\small $\Gamma$}}$_{n-1} (\sigma)$  that 
corresponds to the $reduced$ Gassner representation, is 
 an $(n-1) \times (n-1)$ matrix which is obtained from that 
       abelianized Jacobian  matrix by applying a suitable
conjugation,  
 and deleting the last row and the last column of the form 
$(0, ..., 0, 1)$. 
\smallskip 

Therefore, if we
 denote the free IA-automorphism  corresponding (under 
the $unreduced$ Gassner representation) 
to the braid
$~\sigma~$ by the same letter, the above  condition takes the form
$~rank(J_{\sigma}^{a}  - I)  \le n - 2$,  which points to part
(ii) of our Theorem 1.1. 
\smallskip 

 Thus, we have: 
\smallskip

\noindent {\bf Corollary 1.6.} The Alexander polynomial of a link
$~\hat{\sigma}~$ is zero   if and  only if the free automorphism
which corresponds 
 to the braid $~\sigma$, has a non-trivial fixed point $g \in F'$ 
 modulo $F''$, i.e., $~\sigma(g) = g ~(mod ~F'')$.\\

\noindent {\bf 2. Preliminaries  }
\bigskip

  \indent      Let  $ZF$  be the integral group ring of the 
group  $F$ 
and  $\Delta$   its augmentation ideal, that is, the kernel of the
natural homomorphism  $\epsilon$ : 
$ZF \to Z$.  More generally, when  $R < F$  is a normal 
subgroup  of   $F$,  we denote by $\Delta_R$   the ideal of  $ZF$ 
generated by all elements  of  the  form  
$(r - 1)$,  $r \in R$.  It is the kernel of the natural homomorphism 
$\epsilon_R$ : $ZF \to Z(F/R)$. 
\smallskip 

  \indent     The ideal  $\Delta$  is a free left  $ZF$-module with
a free basis ${\{}(x_i  - 1){\}}$,  $1\le i \le n$, and left Fox
derivations $ d_i $  are projections  to 
the corresponding free cyclic direct summands.  Thus any element 
$ u \in \Delta$ can be uniquely written in the form  
$u =  \sum_{i=1}^n  d_i (u) (x_i - 1)$. 

 \indent        One  can  extend  these  derivations  linearly  to 
 the  whole   $ZF$  by setting $ d_i'(1) = d_i (1) = 0$.  
\smallskip

      The next lemma is an immediate consequence of the definitions.
 
\smallskip

\noindent {\bf Lemma 2.1.}  Let  $J$  be an arbitrary right 
ideal of  $ZF$,   and  let  $u \in \Delta$.  Then  
 $u \in J \Delta$   if and only if $d_i(u) \in J$  for each  $i, 1
\le i \le n$. 

\smallskip

 \indent      Proof of the next lemma can be found in \cite {Fox}. 
\smallskip

\noindent {\bf Lemma 2.2.}  Let  $R$  be a normal subgroup of  $F$,
and let  $y \in F$.  Then   $y - 1 \in$ 

\noindent 
 $\Delta_R \Delta$  if and only if  $y \in R'$. 
\smallskip

 \indent      We also need the ``chain rule" for Fox derivations
(see \cite {Fox}): 
\smallskip 

\noindent {\bf Lemma 2.3.} Let  $\phi$  be  an  endomorphism  of  
$F$   (it  can  be  linearly extended to  $ZF$) defined by 
$\phi(x_k) = y_k , 1 \le k \le n$,  and let  $v =  \varphi(u) $ 
for some  $u, v \in ZF$.  Then: 

\begin{center}
    $ d_j(v) =   \sum_{k=1}^n \varphi(d_k(u))d_j(y_k)$.   
\end{center} 
                   
\smallskip

      Lemma 2.3 implies the  following  product  rule  for  the  
Jacobian  matrices which looks exactly the same as in  the  ``usual" 
 situation  of  analytic functions and Leibnitz derivations: if  
$\varphi$   and   $\psi$   are  two 
 endomorphisms of $F$, then 

\begin{center}
 $J_{\varphi(\psi)} = \psi(J_\varphi) \cdot J_\psi$.
\end{center} 
 
  If furthermore   $\varphi$    and   $\psi$   are
 IA-endomorphisms, then considering abelianization of this  product 
rule yields: 

\begin{center}
 $J_{\varphi(\psi)}^a = J_\varphi^a \cdot J_\psi^a$.
\end{center} 

 In particular,  there  is  a  (faithful)  representation  of  
metabelian  IA-automorphisms by matrices from  $GL_n (Z(F/F'))$ -
cf. \cite {Ba}. 

     Another version  of Lemma 2.3 is:  if  $u, v \in ZF$  and   
 $v = \varphi(u)$,   then 
\begin{eqnarray}
(d_1(v),...,d_n(v))=(\varphi(d_1(u),..., \varphi(d_n(u))) \cdot 
J_\varphi. 
\end{eqnarray} 
     
     Finally, we recall a well-known action via conjugation  of  a  group 
 ring  $Z(F/R)$  on abelian group  $R/R'$  under  which   $R/R'$  
becomes  the 
 (left) relation module of the group  $F/R$.  Namely, if  $g \in
F/R $ and  $r \in  R$,  then  $g$  acts on  $rR'$  by taking it to  
$g r g^{-1}   R'$;   then  this 
 action is extended to the whole group ring  $Z(F/R)$  by 
 linearity.   The  result of the action of  $u \in Z(F/R)$  on 
$rR'$  is denoted by  $r^u R'$.  Then  (see \cite {G}): 

\smallskip 

\noindent {\bf  Lemma 2.4.} For  $r \in R$,  $~u \in Z(F/R)$,  one
has  $d_i(r^u) = u \cdot d_i(r)  ~(mod ~\Delta_R )$. 

\newpage

 \noindent {\bf 3. Proofs }
\bigskip
 
 \noindent {\bf  Proof of Theorem 1.1. }

 \noindent {\bf  (i) }  Let  $\varphi$  be an
 IA-endomorphism of the group  $M$; $~\varphi(x_i) = 
y_i$,  $ ~1 \le i \le n$. 
\smallskip 

  Then, by the    definition of Fox derivatives, we have 
$y_i - 1 =  \sum_{k=1}^n  d_k (y_i) (x_k - 1)~$ for any $~i~$, $ ~1
\le i  \le n$. 

 Since $~\varphi$ is IA,   abelianizing  this equality gives 
 $~x_i^a - 1 =  \sum_{k=1}^n  d_k^a (y_i) (x_k^a - 1),~$ or, 
 in the matrix form:   

\begin{center}
$((x_1^a - 1),...,(x_n^a - 1))^t = J_\varphi^a
\cdot((x_1^a - 1),...,(x_n^a - 1))^t$,  
\end{center} 

 \noindent  where $~^t~$ means taking the 
transpose 
 (i.e., we consider a column, not a row). This is equivalent to 
\begin{center} 
$(J_\varphi^a - I) \cdot ((x_1^a - 1),...,(x_n^a - 1))^t = 0$, 
\end{center} 
hence the columns of the
matrix  $(J_\varphi^a - I)$  are 
 dependent, so  $det(J_\varphi^a - I) = 0$. 

\bigskip

 \noindent {\bf  (ii) }  Let $\varphi$ take $~x_i$ to $x_i s_i$, 
$~s_i \in M'$,  $~1 \le i \le n$.  
For notational convenience, consider the endomorphism $\varphi$
being  lifted to an endomorphism of $F$. We shall use the same
 letter  $\varphi$ for this lifted endomorphism. When $u \in F$, we
denote the abelianization of $u$ (i.e., the image in the group
$A = F/F'$) by $~u^a$.  This agreement also applies to elements of
the group ring $ZF$.  
\smallskip

     If  $det(J_\varphi^a - I) = 0$,  then the rows of the matrix   
$(J_\varphi^a - I)$   are  dependent over the group ring  $ZA = 
Z(F/F')$.  Make a new matrix $J'$ upon multiplying the $~k$th row of 
the matrix   $(J_\varphi^a - I)$ by $~x_k^a$, $~1 \le k \le n$. It is
clear that the rows of $J'$ are  dependent over the group ring  
$ZA$, too. 
\smallskip 

  Now, $J'$ is the abelianized Jacobian 
 matrix of the endomorphism $\varphi'$ that takes  $x_i$ to $s_i$, 
 $~1 \le i \le n$. Indeed, $d_j(x_i s_i) = \delta_{ij} +
x_id_j(s_i)$, so that  $~x_id_j(\varphi'(x_i))  = d_j(\varphi(x_i))
 - \delta_{ij}$. 
\smallskip 

 Thus, the abelianized Jacobian  matrix $J_{\varphi'}^a$ has 
 dependent rows (over the group ring  $ZA$). We 
can write this  dependence as follows: 
\begin{eqnarray}
\sum_{k=1}^n u_k \cdot  d_i^a(s_k) = 0 
\end{eqnarray} 
 \noindent in $ZA$ for some $u_k \in ZA$, not all  of them 0, and $i$ runs
from 1 through $n$. 

By Lemmas 2.1, 2.2, 2.4, the system (2) is equivalent to the following
relation: 
\begin{center}
 $\displaystyle{\prod_{k=1}^n} s_k^{u_k} = 1$ 
\end{center} 
 \noindent in the group $M$. 
\smallskip 

 If $rank ~J_{\varphi'}^a = m \le n - 2$, then there are precisely 
$m$ independent elements (let us say, $s_1, ..., s_m$) among $s_k$
(in other words, these  $s_1, ..., s_m$ generate a free submodule
of the relation module of the group $A$), 
 and for any $j$, $~m + 1 \le j \le n,~$ we have a relation of the 
 form  
 \begin{eqnarray}
s_j^{w_j} = \prod_{k=1}^m s_k^{v_{kj}}
\end{eqnarray} 
 \noindent for some  $~w_j, ~v_{kj} \in ZA$. 
\smallskip 

Now, we are going to find a non-trivial fixed point $g \in M'$ of 
 the endomorphism $\varphi$, in the form 
 \begin{eqnarray}
g = [x_1, x_2]^{z_1}\cdot 
[x_2, x_3]^{z_2} \cdot ... \cdot [ x_{n - 1}, x_n]^{z_{n - 1}} 
\end{eqnarray} 
\noindent  for $z_k \in ZA$  (to be found). 
\smallskip 

 We need a couple of  simple observations:  
\smallskip

\noindent {\bf Lemma 3.1. }{\bf (i)} Let $\varphi(g) = g$ for  some
$g \in M'$. Then  
 $\varphi(g^u) = g^u$ for  any 

\noindent $u \in ZA$. 
\smallskip 

\noindent {\bf (ii)} For any element $h \in M'$, one has 
$h^w = g$  for some $w \in ZA$ and for 
$g \in M'$ of the form (4). 
\smallskip

\noindent {\bf   Proof.}  {\bf (i)}  Since $\varphi$ is IA, we have 
 $\varphi(u) = u$ in $ZA$. The result follows. 
\smallskip 

\noindent {\bf (ii)} It is clear that the elements $[x_i, x_j]$, 
$1  \le i < j \le n,~$ generate $M'$ as a normal subgroup of $M$. 
It is well-known that the elements $[x_1, x_2], ..., [x_{n - 1},
x_n]$ generate a free submodule (of rank $(n - 1)$) 
of the relation module of the group $A$. On the other hand, 
 there are no free submodule of rank $~n~$ in this relation module. 
 This means for any pair $~i, j~$ of indices, one has 
$[x_i, x_j]^{v_{ij}} = [x_1, x_2]^{u_1}\cdot 
[x_2, x_3]^{u_2} \cdot ... \cdot [ x_{n - 1}, x_n]^{u_{n - 1}}~$
for some $~v_{ij}, ~u_k \in ZA$. The result follows. 
 \medskip

 We continue now with part (ii) of Theorem 1.1.
 From $\varphi(g) = g$, we get 
\begin{center}
$[s_1, x_2]^{z_1}\cdot [x_1, s_2]^{z_1}\cdot ... 
\cdot [s_{n - 1}, x_n]^{z_{n - 1}} \cdot  [x_{n - 1}, s_n]^{z_{n -
1}} = 1$, 
\end{center} 
\noindent  or, equivalently: 
 \begin{eqnarray}
s_1^{(x_2 -1)z_1} \cdot s_2^{(1 - x_1)z_1} \cdot ... 
 \cdot s_{n - 1}^{(x_n - 1)z_{n - 1}} \cdot s_n^{(1 - x_{n - 1})
 z_{n -1}} = 1.
\end{eqnarray} 

 Conjugate both sides of (5) by $w = \displaystyle{\prod_{k=m+1}^n}
w_j$,  where  $w_j \in ZA$ come from (3). Then, replace every
$s_j^{w_j}, ~j \ge m+1,~$ in (5), with the corresponding element on
the right-hand side of (3). This gives 
 \begin{eqnarray}
 s_1^{z'_1}\cdot  ... \cdot s_m^{z'_m} = 1
\end{eqnarray} 
\noindent for some $z'_k \in ZA$, each of which is a $ZA$-linear 
combination  of $z_i$, $ ~1 \le i \le n-1$. 
\smallskip 

 Since we have chosen $s_1, ..., s_m$ so that they generate a  free
submodule of the relation module of the group $A$, the equation (6) 
 is equivalent to a system of equations 
 \begin{eqnarray}
 z'_k = 0, ~1 \le k  \le m. 
\end{eqnarray} 

 This is a system of $~m~$ $homogeneous$ 
 $~ZA$-linear equations in $~(n-1) > m$  unknowns $z_1, ... , z_{n
 -1}$. It is well-known that a system like that has a non-trivial solution 
over $ZA$ (since $ZA$ is a commutative domain). This completes the
proof. 
\bigskip

\noindent {\bf (iii)} The ``if" part is obvious. The  ``only if" 
 part follows from  Lemma 3.1. 
\bigskip

 The following example shows how subtle a situation might be in the
case when $~rank(J_{\varphi}^{a}  - I) = n - 1$. 
 \medskip

\noindent {\bf   Proposition 3.2.} In the group $M_2$, 
  let $\varphi$ take $~x_1$ to $x_1 s$, $~x_2$ to  $~x_2 s^{-1}~$ 
for some $s \in M'_2, ~s \ne 1~$. Then $\varphi$ has no non-trivial 
fixed points. 
\smallskip 

\noindent {\bf   Proof.}  First we show that $\varphi$ has no 
 non-trivial  fixed points in $M'_2$. Suppose $\varphi(g) = g; 
 ~g \in M'_2$. Since $~g~$ has the form  $~[x_1, x_2]^z~$ for some 
 $z \in ZA$, we have (cf. (5)): 
 \begin{eqnarray}
 [x_1, x_2]^z = [x_1 s, ~x_2 s^{-1}]^z = [x_1, x_2]^u\cdot s^{(x_2
-1)u + (x_1 -1) z}.
\end{eqnarray} 

 It follows that $s^{(x_2-1)z + (x_1 -1) z} = 1,~$ hence 
 $(x_2-1)\cdot z + (x_1 -1)\cdot z = $

\noindent $(x_2 + x_1 - 2)\cdot z = 0$ 
 in the group ring $ZA$, which is only possible if $z = 0$. Thus, 
 $g = 1$. 
\smallskip 

 Now let $g \in M_2$ be an arbitrary element. Then  $~g~$ has the 
form  $x_1^m x_2^n [x_1, x_2]^z~$ for some  $z \in ZA$. Let us see 
 first what the image of $~x_1^m x_2^n ~$ looks like. 

 We have: $\varphi(x_1^m x_2^n) = x_1 \cdot s\cdot  ... \cdot x_1 \cdot 
s \cdot x_2 \cdot s^{-1}\cdot ... \cdot x_2 \cdot s^{-1} = h$. Now we
apply to this element $~h~$ the following ``collection process":
first we 
 collect all the $x_1$ on the left by permuting them with $~s$. 
 This gives: $h = x_1^m\cdot s^{x_1^{m-1}} \cdot ... \cdot s^{x_1} 
 \cdot s \cdot x_2 \cdot s^{-1}\cdot ...   \cdot x_2 \cdot s^{-1}$.
Then, in the same manner, we collect all the $x_2$ on the right of
$x_1^m$: 

\begin{center}
$h = x_1^m x_2^n \cdot s^{x_1^{m-1} x_2^n}  \cdot ... \cdot s^{x_1 
x_2^n} \cdot s^{x_2^n} \cdot s^{-x_2^{n-1}}  \cdot ... \cdot 
 s^{-x_2} \cdot s^{-1}$, 
\end{center} 

\noindent  or, equivalently: 
 \begin{eqnarray}
h = x_1^m x_2^n \cdot s^{x_1^{m-1} x_2^n +
... + x_1 x_2^n + x_2^n - x_2^{n-1} - ... - x_2 - 1}.  
\end{eqnarray} 
 Now  write down the whole image of $g$:  $~\varphi(g) = 
 \varphi(x_1^m x_2^n \cdot [x_1, x_2]^z) = h \cdot [x_1,
x_2]^z\cdot  s^{(x_2 -1)z + (x_1 -1) z}~$  (cf. (8)). If
$~\varphi(g) = g$, then  combining this with (9) yields: 

\begin{center}
$s^{x_1^{m-1} x_2^n +
... + x_1 x_2^n + x_2^n - x_2^{n-1} - ... - x_2 - 1 + (x_2 -1)z + 
(x_1 -1) z} = 1$, 
\end{center} 

\noindent  or, equivalently: 
 \begin{eqnarray}
x_1^{m-1} x_2^n + ... + x_1 x_2^n + x_2^n - x_2^{n-1} - ... - x_2 -
1 + (x_2 + x_1 - 2)\cdot z  = 0. 
\end{eqnarray} 
 All we have to do now is to show that  $~w = x_1^{m-1} x_2^n + ...
+  x_1 x_2^n + x_2^n - x_2^{n-1} - ... - x_2 - 1~$ is not 
  divisible by $~v = x_2 + x_1 - 2~$  in the ring $ZA$. This is
easy to see upon setting $~x_1 = -1$;  then $~w' = (-1)^m x_2^n - 
 x_2^{n-1} - ... - x_2 - 1;~$  $~v' = x_2 - 3$. Since $~x_2 = 3$ 
 is not a root of the polynomial  $~w'$,  the result follows. This 
 completes the proof of Proposition 3.2. 
\bigskip

\noindent {\bf   Proof of Theorem  1.2.} Let  $\varphi$  be an 
IA-endomorphism of the   group   $M = M_n$. 

First of all, we compute 
the rank of the matrix $J_{\varphi}^{a}~$. If it is not equal to 
 $(n-1)$, then we just refer to Theorem  1.1 (i), (ii). 
\smallskip 

 Suppose $~rank(J_{\varphi}^{a}  - I) = n - 1$. To find out if
there is a non-trivial fixed point of $\varphi$ inside
the commutator subgroup $M'$, we consider a system (7);     but 
 this time, it is  a system of $~(n-1)~$ homogeneous  
 $~ZA$-linear equations in $~(n-1)~$  unknowns $z_1, ... , z_{n
 -1}$. To find out if it has a non-trivial solution (which happens 
 if and  only if $\varphi$ has a non-trivial fixed point inside 
 $M'$), we just compute the corresponding determinant and see if it
is equal to 0. 
\smallskip 

 A somewhat more difficult problem is to find out if
there is a non-trivial fixed point of $\varphi~$ $outside  $ $M'$.
 In this case, we proceed as in the proof of Proposition 3.2, but 
 instead of having just one equation of the form (10), we'll 
 have a system of $~(n-1)~$ $~ZA$-linear equations (they are no
 longer homogeneous!) in $~(n-1)~$  unknowns. 

 Again, since $ZA$ is a commutative domain, we can resolve this 
 system (by using the Kramer's formula). To apply the Kramer's
formula, all we need is to be able to find out for a given pair of 
 polynomials, whether or not one of them is divisible by another.
 But this is equivalent to asking whether or not one of them 
 belongs to the ideal generated by another one. Algorithms like
that do exist (in particular, for $ZA$); they are based on what is
known as Gr\"{o}bner reduction process. 
\smallskip 

 In any case, the existing algorithms in (Laurent) polynomial
algebras not only tell us whether or not a given system of 
 $~ZA$-linear equations has a solution, but if it does, they give
  a solution (although we may not be able to find $all$ of them). 
This  completes the proof of Theorem  1.2.
\bigskip

\noindent {\bf   Proof of Proposition 1.3.}  First of all, we note
that the group $~Fix~\varphi ~\cap ~M'~$ is abelian; therefore, if it
were finitely generated, then its every subgroup would be finitely 
generated, too. 
\smallskip 

 Now suppose $g \in Fix~\varphi ~\cap ~M'$;~$ ~g \ne 1$. Then, by
Lemma 3.1  (i), $g^u \in Fix~\varphi ~\cap ~M'~$ for any $u \in
ZA~$ (i.e., $Fix~\varphi ~\cap ~M'$ is a  normal  subgroup  of 
$M$). We are going to  show that the subgroup of $M'$ generated by
all  $g^{x_1^k},  ~k \in Z$,  is not finitely generated. 

 By means of contradiction, suppose it is generated by  
 $~g^{x_1}, ... , g^{x_1^m}, ~m >0$. 
  But every element from this finitely generated group has a form 
 $g^p$  for some Laurent polynomial $~p = p(x_1),~$ whose degree 
 does not exceed $~m$. Therefore, we don't have the element 
 $g^{x_1^{m+1}}$ in this group,  hence a contradiction. 
 \medskip

\noindent {\bf   Proof of Theorem  1.4.} By means of contradiction, 
suppose $H = Fix ~\varphi $ is a non-cyclic subgroup of  
$M$ generated by  $h_1, ... , h_r,$  $~r < n$. Then there is a
generator of the group $M$, say  $x_1$,  such 
  that $x_1^m \not\in H \cdot M'~$  for any $~m \ge 1$. 
\smallskip 

 By Lemma 3.1 (i), the group $S = H ~\cap ~ M'$ is normal in $M$. If 
$H$ is non-cyclic, then $S$ is non-trivial since it 
contains a non-trivial subgroup $H'~$ (note that $H \not\subseteq
M'~$ since  otherwise, $H$ would be infinitely generated by 
Proposition 1.3). Let $s \in S$, $~s \ne 1$. 
Then for any $~m \ge 1$, we have an equality of the form 
\begin{center}
$\displaystyle{\prod_{j}h_{i_j}^{c_{i_j,m}}} = s^{x_1^m}$
\end{center} 
 for some integers $c_{i_j,m}$. Let $d_k^a(s) \ne 0$;  then (by 
 Lemma 2.4) 
applying the derivation  $d_k^a~$ to both sides of the last 
  equality  gives (in the group ring $ZA$): 
 \begin{eqnarray}
\sum_{q}  \prod_{j}h_{i_j}^{c_{i_j,m}} \cdot c_{i_q,m}\cdot 
  d_k^a(h_{i_q}) = x_1^m \cdot d_k^a(s)  
\end{eqnarray} 
for some collection of indices  $~j, q~$  (of no particular 
   importance to us). 
\smallskip 

 When $~m~$ runs through 1 to $\infty$, we  are
encountering representatives of 
infinitely many  distinct cosets of the group $H\cdot M'/M'~$ (as a
subgroup of $~M/M'$)
 in the supports of elements on the right-hand side of (11). This is
  due 
 to the  condition $x_1^m \not\in H \cdot M'~$  for any $~m \ge 1$
- see above. (By support of a group ring element $~u \in ZG,~$ 
$~u = \sum c_g \cdot g,~$  we mean the set ${\{} g, ~c_g \ne
0{\}}$).  
\smallskip 

 At the same time, the collection of coset representatives of 
$H\cdot M'/M'~$ in the supports of elements on the left-hand side of
(11) does not 
  depend on  $~m,~$ and is therefore finite. This contradiction 
 completes the proof of Theorem  1.4.

 \medskip

\noindent {\bf   Proof of Proposition 1.5.}  
   Consider  an (IA-)automorphism  $\varphi$    of the  group   
$M_n$,   $~ n \ge 3, ~$  
 given by  $\varphi(x_1) = x_1[x_2, x_3, x_1];~$  $~\varphi(x_i) = 
x_i,  ~i \ge  2$.  We are going to prove that the group $~Fix
~\varphi~$ is infinitely generated. 
 \smallskip 

 It is clear that $~Fix~\varphi~$ contains a subgroup of $M_n$ 
generated by $x_2, ..., x_n$. Now we are going to look for (other) 
fixed points in the form 
\begin{center}
$g = x_1^k \cdot [x_1, x_2]^{u_1}\cdot 
[x_2, x_3]^{u_2} \cdot ... \cdot [ x_{n - 1}, x_n]^{u_{n - 1}}
\cdot  x_2^{m_1} \cdot ... \cdot x_n^{m_{n-1}}$.
\end{center} 
\indent  Starting with  $~\varphi(g) = g$  and arguing along  the 
  same lines  as in the proof of Theorem  1.1 (ii) and Proposition
3.2, we finally arrive at 
\begin{center}
$ x_1^k  - 1 + (x_1 - 1)(x_2 -1) \cdot u_1 = 0$, 
\end{center} 
\noindent  which is only possible if $~k = u_1 = 0$. 
 \smallskip 

  It follows that $~x_1^m \not\in Fix~\varphi \cdot M'~$  for any
$~m \ge 1$, 
 and applying the argument from the proof of Theorem  1.4 yields 
     the result. \\

\noindent {\bf Acknowledgement}
\smallskip

\indent I am grateful to A.Krasilnikov and V.Roman'kov 
for useful discussions, 
  and to the referee for helpful comments.

\medskip
\noindent 
 Department of Mathematics, University of California, 
Santa Barbara, CA 93106 
 
\smallskip

\noindent {\it e-mail address\/}: shpil@math.ucsb.edu

\end{document}